\documentclass[a4paper]{article}
  \usepackage[latin1]{inputenc}
  \usepackage[dvips]{graphicx}
  \usepackage{multirow}
  \usepackage{amsmath,amsxtra,amscd,amssymb,latexsym,stmaryrd,theorem}

\newtheorem{defi_aux}{Definition}[section]
\newtheorem{lemm_aux}[defi_aux]{Lemma}
\newtheorem{coro_aux}[defi_aux]{Corollary}
\newtheorem{prop_aux}[defi_aux]{Proposition}
\newtheorem{theo_aux}[defi_aux]{Theorem}
\newtheorem{conj_aux}[defi_aux]{Conjecture}
\newtheorem{rema_aux}[defi_aux]{Remark}

\newenvironment{defi}{%
\medskip\begin{defi_aux}}%
{\end{defi_aux}\medskip}
\newenvironment{lemm}{%
\medskip\begin{lemm_aux}}%
{\end{lemm_aux}\medskip}
{\end{coro_aux}\medskip}
\newenvironment{prop}{%
\medskip\begin{prop_aux}}%
{\end{prop_aux}\medskip}
\newenvironment{theo}{%
\medskip\begin{theo_aux}}%
{\end{theo_aux}\medskip}
{\end{conj_aux}\medskip}
\newenvironment{rema}{%
\medskip\begin{rema_aux}}%
{\end{rema_aux}\medskip}

\newcommand{\preuve}{%
\medskip\noindent{\scshape Proof: }}

\newcommand{\finpreuve}[1]{%
\hspace*{\fill}\rule{6pt}{6pt}\hspace{2pt}{\bf #1}
\bigskip}

\newcommand{\titre}[1]{\maketitle}

\newcommand{\signature}{%
\bigskip
\begin{flushright}
UMPA, \'ENS Lyon\\
46, all\'ee d'Italie\\
69\,364 Lyon cedex 07\\
France\\
{\texttt bkloeckn@umpa.ens-lyon.fr}
\end{flushright}}

\newcommand{\remarque}[1]{\smallskip{\small{\slshape\noindent%
Remark: }#1}\medskip}

\newcommand{\etape}{\sl}

\newcommand{\ensemble}[2]{\left\{#1 ; #2\right\}}

\newcommand{\defini}[1]{{\em #1}\index{#1}}

\newcommand{\abs}[1]{\left\vert #1 \right\vert}
\newcommand{\norme}[1]{\left\Vert #1 \right\Vert}

\newcommand{\impart}{\mathrm{Im}\,}
\newcommand{\repart}{\mathrm{Re}\,}

\newcommand{\id}{\textup{Id}}

\newcommand{\mR}{\ensuremath{{\mathbb R}}}

\newcommand{\mC}{\ensuremath{{\mathbb C}}}

\newcommand{\mS}{\ensuremath{{\mathbb S}}}

\newcommand{\mP}{\ensuremath{{\mathbb P}}}
\newcommand{\mH}{\ensuremath{{\mathbb H}}}

\newcommand{\pour}{\quad}

\newcommand{\comp}{\circ}

\newcommand{\priv}{\setminus}

\newcommand{\diffb}[1]{\ensuremath{{\cal C}^{#1}}}

\newcommand{\agit}{{\bf \cdot}}

\newcommand{\SLDR}
  {\ensuremath{\textrm{\textup{SL}}_2\textup{(}\mR\textup{)}}}
\newcommand{\sldr}
  {\ensuremath{\mathfrak{sl}_2\textup{(}\mR\textup{)}}}
\newcommand{\SLDC}
 {\ensuremath{\textrm{\textup{SL}}_2\textup{(}\mC\textup{)}}}

\newcommand{\SODR}
  {\ensuremath{\textrm{\textup{SO}}_2\textup{(}\mR\textup{)}}}

\newcommand{\RS}{\ensuremath{\overline{\mathbb{C}}}}

\newcommand{\conjugue}{\overline}

\newcommand{\refeq}[1]{\textup{(\ref{#1})}}

\newcommand{\lat}[1]{{\it #1}}

\newcommand{\etoile}{\ensuremath{\null^{\ast}}}

\newcommand{\insiste}[1]{{\sl #1}}

\newcommand{\adherence}[1]{\overline{#1}}

  \title{On differentiable compactifications of the hyperbolic plane 
         and algebraic actions of \SLDR\ on surfaces}
  \author{Beno\^{\i}t Kloeckner}

\begin{document}

%%%%%%%%%%%%%%%%%%%%%%%%%%%%%%%%%%%%%%%%%%%%%%%%%%%%%%%%%%%%%%%%%%%%
%%%%%%%%%%%%%%%%%%%%%%%%%%%%%% Titre %%%%%%%%%%%%%%%%%%%%%%%%%%%%%%%
%%%%%%%%%%%%%%%%%%%%%%%%%%%%%%%%%%%%%%%%%%%%%%%%%%%%%%%%%%%%%%%%%%%%

\titre{On differentiable compactification of the hyperbolic plane
       and algebraic actions of \SLDR\ on surfaces.}

%%%%%%%%%%%%%%%%%%%%%%%%%%%%%%%%%%%%%%%%%%%%%%%%%%%%%%%%%%%%%%%%%%%%%
%%%%%%%%%%%%%%%%%%%%%%%%%%%%%%%%%%%%%%%%%%%%%%%%%%%%%%%%%%%%%%%%%%%%%
\section*{Introduction}
%%%%%%%%%%%%%%%%%%%%%%%%%%%%%%%%%%%%%%%%%%%%%%%%%%%%%%%%%%%%%%%%%%%%%
%%%%%%%%%%%%%%%%%%%%%%%%%%%%%%%%%%%%%%%%%%%%%%%%%%%%%%%%%%%%%%%%%%%%%

An important role played by \SLDR\ is its isometric action on the
hyperbolic plane $\mH^2$, which can be described as the homogeneous 
space $\SLDR/\SODR$, denoted by $\mathcal{E}$.  
This action is real analytic and is, 
up to analytic change of coordinates, the only real analytic 
transitive action of \SLDR\ on the open disk.

The notion of asymptotic geodesics is a means of understanding the
behaviour at infinity of this action, that is to say of giving a
natural topological equivariant compactification of this action
to an action on the closed disk.  
% Moreover, this compactification is topologically unique. 

One can ask whether there is a \insiste{differentiable}
equivariant compactification of this action into the closed disk.
The answer is positive, and there are two well known ways to 
achieve such a compactification.

The restriction to \SLDR\ of the natural action of \SLDC\ on the
Riemann sphere \RS\ has three orbits:
two open hemispheres and between them a great circle.
Considering the union of one open orbit and the circle, one gets 
an analytic equivariant compactification of $\mathcal{E}$. We call
it the \defini{conformal action}.  It corresponds to the 
continuous prolongation to the closed unit 
disk of the \SLDR\ 
action on Poincar\'e's disk.

One can also realize the hyperbolic plane by taking a
lorentzian scalar product $Q$ on $\mR^3$: \SLDR\ acts isometrically on 
$(\mR^3,Q)$, and when one projectivizes $\mR^3$
it gives an analytic action of \SLDR\ on $\mR\mP^2$ with three
orbits: an open disk (which is the hyperbolic plane), 
an open Moebius strip and between them a circle. By taking the action
of \SLDR\ on the union of the disk and the circle we get another
analytic equivariant compactification of $\mathcal{E}$, called
the \defini{projective action}. It corresponds to the continuous
prolongation
to the closed unit disk of the \SLDR\ action on Klein's disk.

By uniqueness, we know that these two compactifications are 
topologically conjugate. However it is easy to check
the following surely known but striking fact:

\begin{prop}
The conformal and projective actions are not
\diffb{1} conjugate, and in particular not $\diffb{\omega}$ conjugate.
\end{prop}

\preuve 
if we choose a point $x$
of the disk boundary and consider in
Poin\-caré's model the closure of the geodesics which have $x$
as an endpoint, we see that all of them are tangent, hence the
differential in $x$ of the conformal action of the
parabolic elements of \SLDR\ which fix $x$ have a common 
proper direction
transversal to the boundary. 

If we now consider the same 
geodesics
in Klein's model, we see that no two of them are tangent and
for each line of the tangent space in $x$, there is a closure 
of a geodesic tangent to it. Hence the differential in
$x$ of the projective action of a parabolic element of
\SLDR\ which fixes $x$ has no proper direction transversal to
the boundary.
\finpreuve

One can ask whether these two compactifications are the only ones. 
The answer, stated in a different way, was given by 
Schneider \cite{Schneider} 
and Stowe \cite{Stowe}:
there exists a countable family of non-equivalent analytic
compactifications of $\mathcal{E}$, which can be described
in terms of infinitesimal generators (see \ref{AnalComp} page
\pageref{AnalComp}). 
These authors also describe all the
analytic actions of \SLDR\ on compact surfaces with or without
boundary and on $\mR^2$.

However these new actions seem less natural than to the two
compactifications we discussed before, which have well known explicit 
integral models. Both
of these models come in a certain sense from the projectivization 
of a linear representation; they will be called \defini{algebraic} 
in the following sense:

\begin{defi}
Let $k$ be a positive integer, possibly $\infty$ or $\omega$.
An action  $\alpha$ of a Lie group $G$ on a manifold possibly with boundary 
$M$ (where $\alpha$, $G$ and $M$ are assumed to be \diffb{k}) is said
to be \diffb{k}-\defini{algebraic}
if there exists a continuous linear representation $\tilde{\rho}$ of 
$G$ on a real finite dimensional vectorial space $V$ and a \diffb{k} 
embedding 
$\Phi: M\longrightarrow\mP(V)$ such that:
\begin{itemize}
  \item $\Phi(M)$ is a union of orbits for the action $\rho$ induced 
        by $\tilde{\rho}$ on $\mP(V)$,
  \item $\alpha$ coincides with $\rho$ via $\Phi$, that is:
    $$ \Phi\comp\alpha(g) = \rho(g)\comp\Phi 
       \pour\forall g\in G. $$
\end{itemize}
The pair $(\tilde{\rho},\Phi)$ is called a \diffb{k} 
\defini{algebraic realization} of $\alpha$.
\end{defi}

It is obvious that the projective action is algebraic. The
Riemann sphere can be seen as a submanifold of the space of the
2-plans of $\mR^4$ which, as a Grassmanian, can be embedded in a 
real projective space such that the conformal action of \SLDR\ 
extends to the projectivization of a linear representation. So the
conformal action is algebraic too.

By studying the topology of all the algebraic continuous actions of 
\SLDR\ on surfaces and thus determining the regularity of the gluing 
of the orbits we prove (for a precise definition of 
``compactification'' see \ref{AnalComp}):

\begin{theo}\label{TheoEllip}
The conformal and projective actions are the only \diffb{\omega}
compactifications of $\mathcal{E}$ which are algebraic.
\end{theo}

With this material, we are also able to study all the analytic algebraic
actions of \SLDR\ on surfaces and prove:

\begin{theo}
The analytic algebraic actions of \SLDR\ on surfaces (with or without
boundary) consist exactly
of:
\begin{itemize}
  \item the projective action (on $\mR\mP^2$),
  \item the conformal action (on $\mS^2$),
  \item the standard product action on $\mR\mP^1\times\mR\mP^1$,
  \item one action on the projective plane with an open dense orbit,
  \item a countable family of actions on the Klein bottle,
  \item a countable family of actions on the torus with 
        two open cylindric orbits and two circular orbits,
  \item a countable family of actions on the torus with
        four open cylindric orbits and four circular orbits,
\end{itemize}
and of any subaction (\lat{i.e.} union of orbits) of any one of
these actions.
\end{theo}

\remarque{The realization of these actions as algebraic actions 
gives explicit global models for all of them.}

%%%%%%%%%%%%%%%%%%%%%%%%%%%%%%%%%%%%%%%%%%%%%%%%%%%%%%%%%%%%%%%%%
%%%%%%%%%%%%%%%%%%%%%%%%%%%%%%%%%%%%%%%%%%%%%%%%%%%%%%%%%%%%%%%%%
\section{The topology of low dimensional algebraic orbits}
%%%%%%%%%%%%%%%%%%%%%%%%%%%%%%%%%%%%%%%%%%%%%%%%%%%%%%%%%%%%%%%%%
%%%%%%%%%%%%%%%%%%%%%%%%%%%%%%%%%%%%%%%%%%%%%%%%%%%%%%%%%%%%%%%%%

	Our goal is in this section to describe the topology of 
all orbits of dimension less or equal to 2 which appear in the 
projectivization of a finite dimensional linear representation
of \SLDR.

%%%%%%%%%%%%%%%%%%%%%%%%%%%%%%%%%%%%%%%%%%%%%%%%%%%%%%%%%%%%%%%%%
\subsection{Irreducible representations}
%%%%%%%%%%%%%%%%%%%%%%%%%%%%%%%%%%%%%%%%%%%%%%%%%%%%%%%%%%%%%%%%%

All the irreducible representations of \SLDR\ are known; for a 
proof of the following theorem, see \cite{Serre}.

We define a family of linear representations of \SLDR. For 
each non-negative integer $n$, 
$\tilde\rho_n : \SLDR\longrightarrow\mR_n[X,Y]$, where
$\mR_n[X,Y]$ is the vector space of all homogenous polynomials
of degree $n$ in $X$ and $Y$, is given by
$$\tilde\rho_n \begin{pmatrix}a&b\\c&d\end{pmatrix} \agit P(X,Y)
  = P(aX+cY,bX+dY).$$

\begin{theo}
The representation (of dimension $n+1$) $\tilde\rho_n$ is 
irreducible for any non-negative $n$ and any finite-dimensional
irreducible representation of \SLDR\ is of this form.
\end{theo}

%%%%%%%%%%%%%%%%%%%%%%%%%%%%%%%%%%%%%%%%%%%%%%%%%%%%%%%%%%%%%%%%%%
\subsection{Irreducible case}
%%%%%%%%%%%%%%%%%%%%%%%%%%%%%%%%%%%%%%%%%%%%%%%%%%%%%%%%%%%%%%%%%%

We start the study by the irreducible case. 

The irreducible
representation of dimension 1, $\tilde\rho_0$, is trivial:
its associated projective action has one single (fixed !) point. 

The irreducible
representation of dimension 2, $\tilde\rho_1$, gives the obvious
action of \SLDR\ on $\mR\mP^1$, which is transitive. 

The irreducible
representation of dimension 3, $\tilde\rho_2$, gives the projective
action on $\mR\mP^2$, which has three orbits : one open disc, one 
circle and one Moebius strip. We can determine in which orbit lies the
vector line given by a polynomial $P=aX^2+bXY+cY^2$ (we denote such
a line by $[aX^2+bXY+cY^2]$) just by computing the discriminant
$\Delta=b^2-4ac$ (which plays the role of the Lorentzian scalar 
product in the description of the projective action given in the
introduction). The open disk consists of the elements which are not
factorizable over \mR\ (\lat{i.e.} of non-positive 
discriminant). The Moebius strip consists of those which are
factorizable with two distinct factors (\lat{i.e.} of non-negative
discriminant). The circle consists of those which are squares (\lat{i.e.}
of zero discriminant).

We denote by $\mH^+$ the upper half plane in \mC\ and by 
$\partial\mH^+$ its boundary (in Riemann's sphere \RS). We have
a canonical identification between $\partial\mH^+$ and $\mR\mP^1$,
which allows us to identify them.

It is important to notice that, since the map: 
\begin{eqnarray}
\mH^+\sqcup\partial\mH^+ & \longrightarrow & \mP(\mR_3[X,Y])
                                                     \nonumber\\
z                        & \longmapsto     & [(zX+Y)(\conjugue{z}X+Y)]
                                                     \nonumber
\end{eqnarray}
is not differentiable on the boundary, it
is not an analytic parametrization of the closed disk 
(union of the open disk orbit and of the circular orbit) and there
is no reason to think that the conformal and projective actions on the
closed disk are equal up to analytic coordinate change (we already
saw that they are not).

Now we generalize this method for all irreducible representations. We
shall fix a non-negative integer $n$. An element of $\mP(\mR_n[X,Y])$ 
factorizes into the following form:
\begin{equation}\label{factorization}
  \left[\prod_{i=1}^k (t_iX+Y)^{\alpha_i}\prod_{j=1}^l%
  (z_jX+Y)^{\beta_j}(\conjugue{z_j}X+Y)^{\beta_j}\right]
\end{equation}
where $t_i$'s are distinct elements of $\partial\mH^+$,
$z_j$'s are disctinct elements of $\mH^+$ and 
$\sum\alpha_i+2\sum\beta_j=n$.

Note that $t_i$'s are possibly infinite : for example $[\infty X+Y]$ 
denotes the projective element $[X]$.

The form \refeq{factorization} is efficient: since we have
$$ \begin{pmatrix}a&b\\c&d\end{pmatrix}\agit(zX+Y)=
   (cz+d)\left(\frac{az+b}{cz+d}X+Y\right)$$
where $z\in\RS$, the conformal action allows one to
study  all algebraic actions
of \SLDR \insiste{topologically}.

We shall first determine which orbits are of dimension 2 or less.

\begin{lemm}\label{lemme_kl}
The orbit of an element $P$ written under the form
\refeq{factorization} is of dimension 2 or less if 
and only if: $k+2l\leqslant 2$.
\end{lemm}

\preuve
We consider the different cases one by one. By ``isometry''
we shall always mean ``orientation-preserving isometry''.

If $l=1$ and $k=1$, we can write 
$P=[(tX+Y)^{\alpha}(zX+Y)^{\beta}(\conjugue{z}X+Y)^{\beta}]$
and the stabilizer of P is the set of the isometries
of $\mH^+$ (with the hyperbolic metric) which fix the point $z$ 
and the point of the boundary $t$, and hence consist only of the 
identity \id.
Thus the orbit of $P$ is of the same dimension than \SLDR, \lat{i.e.}
3.

If $l\geqslant 1$ and $k\geqslant 1$, the same conclusion holds.

If $l\geqslant2$, an element of the component of \id\ in the 
stabilizer of $P$ must fix at least two points of $\mH^+$,
hence it is discrete and the orbit of $P$ is of dimension 3.

If $k\geqslant3$, an element of the component of \id\ of the 
stabilizer of $P$ must fix at least three points of the
boundary $\partial\mH^+$, hence the same conclusion holds.

If $l=0$ and $k=1$, the stabilizer of $P$ is the set of the 
isometries of $\mH^+$ which fix one given point (the only root of a
representative polynomial for $P$) of the boundary, 
hence its dimension is 2. Thus the orbit of $P$ is one-dimensional.

If $l=0$ and $k=2$, the stabilizer of $P$ is the set of the
isometries of $\mH^+$ which fix two given points of the boundary, 
hence it is one-dimensional. Thus the dimension of
the orbit of $P$ is 2.

If $l=1$ and $k=0$ the stabilizer of $P$ is the set of the
isometries of $\mH^+$ which fix one given point, 
hence it is one-dimensional. Thus the dimension of the orbit 
of $P$ is 2.
\finpreuve

We have three cases of low dimensional orbits, namely the 
\defini{elliptic} case 
($l=1$ and $k=0$), the \defini{parabolic} case ($l=0$ and $k=1$) 
and the \defini{hyperbolic} case ($l=0$ and $k=2$).

\begin{prop}
The topology of an orbit of dimension 2 or less of the action
$\rho_n$ (obtained by projectivizing $\tilde\rho_n$) is 
given by the factorized form \refeq{factorization} of any one
of its elements $P$ in the following way:
\begin{enumerate}
  \item \label{cas_pc} if $l=0$ and $k=1$: the orbit of $P$ 
        is a circle
        $$ \ensemble{[(tX+Y)^n]}{t\in\partial\mH^+}. $$
        There is only one such orbit,
  \item \label{cas_hm} if $l=0$, $k=2$ and $\alpha_1=\alpha_2$: 
        the orbit of $P$ is a Moebius strip
        $$ \ensemble{[(t_1X+Y)^{\alpha}
           (t_2X+Y)^{\alpha}]}{t_1\neq t_2\in\partial\mH^+} $$
        where $t_1$ and $t_2$ play the same role.
        There is one such orbit if $n$ is even, none if $n$
        is odd,
  \item \label{cas_hc} if $l=0$, $k=2$ and $\alpha_1\neq\alpha_2$: 
        the orbit of $P$ is a cylinder
        $$ \ensemble{[(t_1X+Y)^{\alpha_1}
           (t_2X+Y)^{\alpha_2}]}{t_1\neq t_2\in\partial\mH^+} $$
        where $t_1$ and $t_2$ play non-symmetric roles (inverting them
        maps an element of the orbit to another). There are
        $\frac{n-1}{2}$ such orbits if $n$ is odd,
        $\frac{n-2}{2}$ if $n$ is even,
  \item \label{cas_e} if $l=1$ and $k=0$: the orbit of $P$ 
        is a disc
        $$ \ensemble{[(zX+Y)^{\beta}
           (\conjugue{z}X+Y)^{\beta}]}{z\in\mH^+}. $$
        There is one such orbit if $n$ is even, none if
        $n$ is odd.
\end{enumerate}
\end{prop}

\preuve
As \SLDR\ is transitive on $\mH^+$ and doubly transitive on 
$\partial\mH^+$,
each set described here is an orbit. Thanks to Lemma \ref{lemme_kl}
there is no other case than the four mentionned. The computation of the 
number of orbits is easy with the condition 
$\sum\alpha_i+2\sum\beta_j=n$.

All we have to prove is that the topology of each of these sets is as 
claimed. The cases \ref{cas_pc}, \ref{cas_hm}, \ref{cas_e}
can be deduced from the study of $\rho_2$ since the map
\begin{eqnarray}
\mP(\mR_m[X,Y]) & \longrightarrow & \mP(\mR_{\alpha m}[X,Y])
                                                      \nonumber \\
\null[P]            & \longmapsto     & [P^{\alpha}]
                                                       \nonumber
\end{eqnarray}
is a homeomorphism on its image.

The case \ref{cas_hc} reduces to the elementary fact that
$$\ensemble{(x,y)\in\mS^1\times\mS^1}{x\neq y}$$
is a cylinder.
\finpreuve

%%%%%%%%%%%%%%%%%%%%%%%%%%%%%%%%%%%%%%%%%%%%%%%%%%%%%%%%%%%%%%%%%%%%%%%%%
\subsection{Notations for the reducible case}
%%%%%%%%%%%%%%%%%%%%%%%%%%%%%%%%%%%%%%%%%%%%%%%%%%%%%%%%%%%%%%%%%%%%%%%%%

We shall now consider the reducible representations of \SLDR. Since it is a 
semi-simple Lie group, its finite-dimensional representations are 
sums of irreducible representations. If we consider a representation
$\tilde\rho$, we can write: 
$\tilde\rho=\tilde\rho_{n_1}\oplus\tilde\rho_{n_2}\oplus\dots
\tilde\rho_{n_p}$ for some $n_1,\dots,n_p$.

We denote by
$V=\mR_{n_1}[X,Y]\oplus\mR_{n_2}[X,Y]\oplus\dots\oplus\mR_{n_p}[X,Y]$ 
the vector space of $\tilde\rho$. Up to a permutation, we can assume 
that $n_1\geqslant n_2\geqslant\dots\geqslant n_p$.

Moreover, as we want to consider together all the copies of a given
irreducible representation which appears in $\tilde\rho$ we set
$I_1=\llbracket i_1=1,i_2-1\rrbracket$, 
$I_2=\llbracket i_2,i_3-1\rrbracket$, $\dots$,
$I_r=\llbracket i_r,i_{r+1}-1=p\rrbracket$ the integer intervals such that:
$$\underbrace{n_1=\dots=n_{i_2-1}}_{I_1}>
  \underbrace{n_{i_2}=\dots=n_{i_3-1}}_{I_2}>\dots>
  \underbrace{n_{i_r}=\dots=n_p}_{I_r}.$$

We say that $I_s$ is \defini{even}, respectively \defini{odd}
if $n_{i_s}$ is even, respectively odd.

We write an element $x$ of $\mP(V)$ under the factorized form:
\begin{equation}\label{factorization_reducible}
x=\left[ u_q\prod_{i=1}^{k_q}(t_q^iX+Y)^{\alpha_q^i}
            \prod_{j=1}^{l_q}(z_q^jX+Y)^{\beta_q^j}
                             (\conjugue{z_q^j}X+Y)^{\beta_q^j}
  \right]_{1\leqslant q\leqslant p}
\end{equation}
where the $u_q$'s are real numbers and for each $q$:
$\sum\alpha_q^i+2\sum\beta_q^j=n_q$. 

We call \defini{support} of $x$ (or of the projective element
$[u_1,\dots,u_p]$) and denote by $I(x)$ the set of all the intervals 
$I_s$ such that there is at least one index $i\in I_s$, 
$u_i\neq 0$. We write $q\in I(x)$ instead of 
$q\in\bigcup_{I\in I(x)} I$.

We say that a support is \defini{even}, respectively \defini{odd}
if all of its elements are even, respectively odd. We define an 
\defini{odd} support the same way.

We denote by $I_+(x)$ the element of the support of $x$ which carries
the greatest dimension (\lat{i.e.} the lowest indices),
$I_-(x)$ the one which carries the lowest dimension. We denote by
$q_+(x)$ (respectively $q_-(x)$) the smallest (respectively
the greatest) index $q$ such that $u_q\neq 0$. We have
$q_+(x)\in I_+(x)$ and $q_-(x)\in I_-(x)$.

When there is no ambiguity, we write $I_+$, $I_-$, $q_+$ and $q_-$
instead of $I_+(x)$, $I_-(x)$, $q_+(x)$ and $q_-(x)$.

We denote by $k(x)$ (or $k$) the number of different $t_q^i$'s of 
$\partial\mH^+$ which arise
in the factorized form \refeq{factorization_reducible} of $x$, and
$l(x)$ (or $l$) the number of different $z_q^j$'s of $\mH^+$.

With these notations we can now generalise the results of the
previous section to reducible representations.

\begin{lemm}
Let $x$ be a element of the projective space $\mP(V)$ whose orbit
is of dimension 2 or less. Then $k(x)+2l(x)\leqslant 2$.
\end{lemm}

\preuve
An element of the identity component of the stabilizer of $x$
is an isometry of $\mH^+$ stabilizing $l(x)$ points 
and $k(x)$ points of the boundary, so we can conclude using
the discussion in the proof of Lemma \ref{lemme_kl}.
\finpreuve

Until the end of the paper, we shall assume there is at least one
index $i$ such that $n_i>1$ (otherwise the action of \SLDR\ is
trivial).

%%%%%%%%%%%%%%%%%%%%%%%%%%%%%%%%%%%%%%%%%%%%%%%%%%%%%%%%%%%%%%%%%%%%%%%%
\subsection{Reducible elliptic case}
%%%%%%%%%%%%%%%%%%%%%%%%%%%%%%%%%%%%%%%%%%%%%%%%%%%%%%%%%%%%%%%%%%%%%%%%

We assume here that $k=0$ and $l=1$, that is to say we consider 
the orbit of an element
$$ x = \left[u_q(zX+Y)^{\frac{n_q}{2}}
       (\conjugue{z}X+Y)^{\frac{n_q}{2}}
       \right]_{1\leqslant q\leqslant p}$$
which must be of even support.

\begin{lemm}
The orbit of an elliptic element is homeomorphic to a disk.
\end{lemm}

\preuve
composing with an element of \SLDR, we can assume $z=\imath$.
Thus the elements of the stabilizer of $x$ are exactly
the matrices
{\footnotesize$\begin{pmatrix} a & b \\ -b & a \end{pmatrix}$}
where $a^2+b^2=1$.

Hence we can parametrize the orbit of $x$ by $z\in\mH^+$.
\finpreuve

%%%%%%%%%%%%%%%%%%%%%%%%%%%%%%%%%%%%%%%%%%%%%%%%%%%%%%%%%%%%%%%%%%%%%%%
\subsection{Reducible parabolic case}\label{Rpc}
%%%%%%%%%%%%%%%%%%%%%%%%%%%%%%%%%%%%%%%%%%%%%%%%%%%%%%%%%%%%%%%%%%%%%%%

Now we shall assume $k=1$ and $l=0$ and consider an element
$x = \left[u_qY^{n_q}\right]$
(after possible composition with an element of \SLDR).

\begin{lemm}
The orbit of a parabolic element with support reduced to a single 
element is homeomorphic to a circle. 

The orbit of a parbolic element with support containing at least
two elements is homeomorphic to a cylinder. 
\end{lemm}

\preuve
if $A =${\footnotesize $\begin{pmatrix} a&b\\c&d \end{pmatrix}$}
stabilizes $x$, thus it stabilizes $0$ when acting projectively
on $\mR\mP^1$ hence $b=0$ (and $d=a^{-1}$).

Moreover we have
$$\begin{pmatrix}a&0\\c&a^{-1}\end{pmatrix}\agit x
  = \left[u_q a^{-n_q}Y^{n_q}\right]_q.$$

{\etape If the support of $x$ consists of one single interval $I_s$}
the condition $b=0$ is sufficient for $A$ to stabilize $x$.
If $d\neq0$,
$$\begin{pmatrix}a&b\\c&d\end{pmatrix}\agit x
  = \left[u_q \left(\frac{b}{d}X+Y\right)^{n_q}\right]_{q\in I_s}$$
else
$$\begin{pmatrix}a&b\\c&0\end{pmatrix}\agit x
  = \left[u_q X^{n_q}\right]_{q\in I_s}$$
Hence the orbit of $x$ is homeomorphic to $\mR\mP^1$.

{\etape If the support of $x$ consists of at least to intervals}
the stabilizer of $x$ consist of the matrices of the form
$A =${\footnotesize $\begin{pmatrix} 1&0\\c&1 \end{pmatrix}$}
hence the orbit is of dimension 2.

If $d\neq0$,
$$\begin{pmatrix}a&b\\c&d\end{pmatrix}\agit x
  = \left[u_q d^{n_q}\left(\frac{b}{d}X+Y\right)^{n_q}\right]_q$$
else
$$\begin{pmatrix}a&b\\c&0\end{pmatrix}\agit x
  = \left[u_q b^{n_q}X^{n_q}\right]_q$$
hence a point of the orbit of $x$ is determined by
$\frac{b}{d}\in\mR\mP^1$ and a real non-zero parameter,
$b$ or $d$. The case $d\neq0$ gives a pair of disjoint copies of
$\mR\times\mR\etoile$ which are glued along $d=0$ into a
cylinder. If the support of $x$ is neither even nor odd
this cylinder is naturally homeomorphic to the orbit of $x$,
otherwise {\footnotesize $\begin{pmatrix}-a&-b\\-c&-d\end{pmatrix}$}%
$\agit x =$ {\footnotesize $\begin{pmatrix}a&b\\c&d\end{pmatrix}$}
and it is naturally a 2-folded covering of the orbit of $x$
which is a cylinder too.
\finpreuve

%%%%%%%%%%%%%%%%%%%%%%%%%%%%%%%%%%%%%%%%%%%%%%%%%%%%%%%%%%%%%%%%%%%
\subsection{Reducible hyperbolic case}\label{Rhc}
%%%%%%%%%%%%%%%%%%%%%%%%%%%%%%%%%%%%%%%%%%%%%%%%%%%%%%%%%%%%%%%%%%%

We shall assume $k=2$ and $l=0$ and consider an element
$$ x=\left[u_qX^{\alpha_q}Y^{n_q-\alpha_q}\right]_q$$
(note that we define $\alpha_q$ only when $u_q\neq 0$).

\begin{lemm}
With the notations of this section, a hyperbolic element has a
2 dimensional orbit if and only if $2\alpha_q-n_q$ is constant,
noted $\delta$. When this condition is satisfied, the orbit is
a Moebius strip if $\delta=0$ and $\alpha_{q_+}-\alpha_q$ is 
even for each $q$, a cylinder otherwise.
\end{lemm}

\preuve
a stabilizing element of $x$ must stabilize $0$ and $\infty$ in
\RS\ hence can be written 
{\footnotesize$\begin{pmatrix}a&0\\0&a^{-1}\end{pmatrix}$}.
As
{\footnotesize$\begin{pmatrix}a&0\\0&a^{-1}\end{pmatrix}$}%
$\agit x = \left[u_q a^{2\alpha_q-n_q}X^{\alpha_q}Y^{n_q-\alpha_q}\right]$
we see that if there are $q_1$, $q_2$
such that $2\alpha_{q_1}-n_{q_1}\neq 2\alpha_{q_1}-n_{q_1}$ thus
the orbit of $x$ is 3-dimensional, and is 2-dimensional otherwise.

We shall assume we are in the latter case.

Thus the image of $x$ under the action of an element $A\in\SLDR$ is
given by the images $t_1$ and $t_2$ of $0$ and $\infty$ under the 
action of $A$ on $\mR\mP^1$. If $\alpha_q=\frac{n_q}{2}$
for all $q$ ($x$ is therefore of even 
support) and $\alpha_{q_+}-\alpha_q$ is even for all $q$ thus exchanging
$t_1$ and $t_2$ gives the same point of the orbit, else it does not.
\finpreuve

%%%%%%%%%%%%%%%%%%%%%%%%%%%%%%%%%%%%%%%%%%%%%%%%%%%%%%%%%%%%%%%%%%%%%%%%
%%%%%%%%%%%%%%%%%%%%%%%%%%%%%%%%%%%%%%%%%%%%%%%%%%%%%%%%%%%%%%%%%%%%%%%%
\section{Closure of low dimensional algebraic orbits}
%%%%%%%%%%%%%%%%%%%%%%%%%%%%%%%%%%%%%%%%%%%%%%%%%%%%%%%%%%%%%%%%%%%%%%%%
%%%%%%%%%%%%%%%%%%%%%%%%%%%%%%%%%%%%%%%%%%%%%%%%%%%%%%%%%%%%%%%%%%%%%%%%

We shall now determine the closures of the orbits. 

By the \defini{border} of an orbit $O$ we mean the set
$\adherence{O}\priv O$.

%%%%%%%%%%%%%%%%%%%%%%%%%%%%%%%%%%%%%%%%%%%%%%%%%%%%%%%%%%%%%%%%%%%%%%%%
\subsection{Elliptic case}
%%%%%%%%%%%%%%%%%%%%%%%%%%%%%%%%%%%%%%%%%%%%%%%%%%%%%%%%%%%%%%%%%%%%%%%%

We shall consider the orbit of the element $x$ which is elliptic, associated
to $\imath$ and $[u_q]_q$, that is :
$ x = \left[u_q (\imath X+Y)^{\frac{n_q}{2}}
                (-\imath X+Y)^{\frac{n_q}{2}} \right]_q.$
\begin{lemm}
The border of the orbit of an elliptic element $x$ associated to a 
projective point $[u_q]_q$ is the circular parabolic orbit of
$[u_q Y^{n_q}]_{q\in I_+(x)}$. The union of these two orbits
is a closed disk.
\end{lemm}

\preuve
we have
$$\begin{pmatrix} a&b\\c&d \end{pmatrix} \agit x =
  \left[u_q\abs{c\imath+d}^{n_q-n_{q_+}}
  \left(\frac{a\imath+b}{c\imath+d}X+Y\right)^{\frac{n_q}{2}}
  \left(\conjugue{\frac{a\imath+b}{c\imath+d}}X
  +Y\right)^{\frac{n_q}{2}}\right]_q$$

Since $ad-bc=1$ we can write:
$$\frac{a\imath+b}{c\imath+d}=
  \frac{ac+bd}{\abs{c\imath+d}^2}
  +\imath \frac{1}{\abs{c\imath+d}^2}$$
thus $\abs{c\imath+d}^2=(\impart z)^{-1}$, and hence the orbit is the set
of the elements
$$x(z) = \left[u_q(\impart z)^{\frac{n_{q_+}-n_q}{2}}
  \left(zX+Y\right)^{\frac{n_q}{2}}
  \left(\conjugue{z}X
  +Y\right)^{\frac{n_q}{2}}\right]_q$$
where $z\in\mH^+$.

If a sequence $(x(z_i))_i$ has a limit in $\mP(V)$, necessarily $(z_i)_i$
has a limit in the closure of $\mH^+$ in \RS. If this limit is in
$\mH^+$ we get a point of the orbit of $x$, otherwise it is a point 
$t\in\partial\mH^+$. In the latter case, if $t$ is finite, $\impart{z_i}$
has limit zero and $(x(z_i))_i$ has limit $[u_q(tX+Y)^{n_q}]_{q\in I_+(x)}$.
If $t=\infty$, $(x(z_i))_i$ has limit $[u_q X^{n_q}]_{q\in I_+(x)}$,
which we can write $[u_q(\infty X+Y)^{n_q}]_{q\in I_+(x)}$.
\finpreuve

%%%%%%%%%%%%%%%%%%%%%%%%%%%%%%%%%%%%%%%%%%%%%%%%%%%%%%%%%%%%%%%%%%%%%%%%%%%%
\subsection{Parabolic case}
%%%%%%%%%%%%%%%%%%%%%%%%%%%%%%%%%%%%%%%%%%%%%%%%%%%%%%%%%%%%%%%%%%%%%%%%%%%%

The circular orbits are closed, so we consider only the two types of
cylindric orbits; as the technic is the same than in the elliptic 
case, we shall not give much detail.

\begin{lemm}
Let $x=[u_q Y^{n_q}]_q$ be of even non-reduced to a single element support.
The border of the cylindric orbit of $x$ is the disjoint union of the
orbits of $[u_q Y^{n_q}]_{q\in I_+(x)}$ and $[u_q Y^{n_q}]_{q\in I_-(x)}$.

If the support of $x$ has a parity (\lat{i.e.} is even or odd), 
the closure of the orbit of $x$ is a 
closed cylinder if $n_{q_-}>0$ and a closed disk if $n_{q_-}=0$.

If the support of $x$ is neither odd nor even, the closure of the orbit
of $x$ is a Klein bottle if  $n_{q_-}>0$ and a projective plane if
$n_{q_-}=0$.
\end{lemm}

\preuve
we shall consider the orbit of an element $x=[u_q Y^{n_q}]_q$ whose support
is even and has at least two elements. This orbit is described in
Section \ref{Rpc}, we can write it under the form:
\begin{eqnarray}
\begin{pmatrix}a&b\\c&d\end{pmatrix}\agit x
    & = &
    \left[u_q d^{n_q-n_{q_\pm}}\left(\frac{b}{d}X+Y\right)^{n_q}\right]_q
    \mbox{ if }d\neq 0,
    \nonumber\\
\begin{pmatrix}a&b\\c&0\end{pmatrix}\agit x
    & = &
    \left[u_q b^{n_q-n_{q_\pm}}X^{n_q}\right]_q
    \nonumber
\end{eqnarray}
where we choose $\pm$ to be $+$ (respectively $-$) 
if we want to study great (respectively
small) values of the real parameter given for a choosen 
$t=\frac{b}{d}\in\mR\mP^1$ by $d$ (or $b$ if $t=\infty$).

For great values, we find a point of the circular orbit of 
$[u_q Y^{n_q}]_{q\in I_+(x)}$, for small ones a point of the
orbit of $[u_q Y^{n_q}]_{q\in I_-(x)}$ (which is a circle
if $n_{q_-(x)}>0$, a single point otherwise).

The way the cylindric orbit is glued on the circles of its border
depends of the parity of the support of $x$: if it has a parity
(\lat{i.e.} is even or odd) the couples $(b,d)$ and $(-b,-d)$ 
of parameters give the same point,
else they give two different points such that if one of them is close to
a point of the border, the other is close to this point too: hence
the cylinder will glue twice on each circle in its border.
\finpreuve

%%%%%%%%%%%%%%%%%%%%%%%%%%%%%%%%%%%%%%%%%%%%%%%%%%%%%%%%%%%%%%%%%%%%%%%%%%
\subsection{Hyperbolic case}
%%%%%%%%%%%%%%%%%%%%%%%%%%%%%%%%%%%%%%%%%%%%%%%%%%%%%%%%%%%%%%%%%%%%%%%%%%

\begin{lemm}
The border of the orbit $O$ of an element
$x=\left[u_qX^{\alpha_q}Y^{n_q-\alpha_q}\right]_q$
(where $2\alpha_q-n_q$ does not depend upon $q$) is the circular orbit
of $[u_q Y^{n_q}]_{q\in I_+(x)}$.

If $O$ is a Moebius strip, its closure is a closed Moebius strip.

If $O$ is a cylinder, its closure is a torus.
\end{lemm}

\preuve
we can write this
orbit as the set of all elements of the form
$$\left[u_q (t_1-t_2)^{\alpha_{q_+}-\alpha_q}
          \left(t_1 X+Y\right)^{\alpha_q}
          \left(t_2 X+Y\right)^{\beta_q}\right]_q$$
$$
  =\left[u_q
          \left(\frac{1}{t_2}-\frac{1}{t_1}\right)^{\alpha_{q_+}-\alpha_q}
          \left(X+\frac{1}{t_1}Y\right)^{\alpha_q}
          \left(X+\frac{1}{t_2}Y\right)^{\beta_q}\right]_q$$
with $t_1,t_2\in\mR\mP^1$. As before, this enables the description
of the border of this orbit.
\finpreuve

%%%%%%%%%%%%%%%%%%%%%%%%%%%%%%%%%%%%%%%%%%%%%%%%%%%%%%%%%%%%%%%%%%%%%%%%%
%%%%%%%%%%%%%%%%%%%%%%%%%%%%%%%%%%%%%%%%%%%%%%%%%%%%%%%%%%%%%%%%%%%%%%%%%
\section{Classification of analytic algebraic action of \SLDR\ on
         surfaces}
%%%%%%%%%%%%%%%%%%%%%%%%%%%%%%%%%%%%%%%%%%%%%%%%%%%%%%%%%%%%%%%%%%%%%%%%%
%%%%%%%%%%%%%%%%%%%%%%%%%%%%%%%%%%%%%%%%%%%%%%%%%%%%%%%%%%%%%%%%%%%%%%%%%

We shall now study the analyticity of the different topological surfaces
obtained as a union of orbits and which are analytically conjugate
(\lat{i.e.} are equal up to an analytic change of coordinates).

%%%%%%%%%%%%%%%%%%%%%%%%%%%%%%%%%%%%%%%%%%%%%%%%%%%%%%%%%%%%%%%%%%%%%%%%%
\subsection{Smoothness of polynomial-parametrized surfaces}
%%%%%%%%%%%%%%%%%%%%%%%%%%%%%%%%%%%%%%%%%%%%%%%%%%%%%%%%%%%%%%%%%%%%%%%%%

We shall use many times the following result, which can be generalized (but
we present here only the 2-dimensional version for simplicity).

\begin{prop}\label{raccourci}
Let 
$P : (x_1,x_2) \longmapsto \left(P_1(x_1,x_2),\dots,P_n(x_1,x_2)\right)$
be a map defined on a neigborhood of 0 in $\mR^2$ where the $P_i$'s are
homogeneous non-constant polynomials. We assume
$P_1$ to be of minimal degree and 
$P_2\notin\mR[P_1]$ of minimal degree among $P_i$'s with that property.
If there exists some $P_i\notin\mR[P_1,P_2]$ then the image $E$ of $P$ is
not a smooth 2-dimensional submanifold of $\mR^n$ (more precisely, $P$
is singular at 0).
\end{prop}

\preuve
Assume that $E$ is a smooth 2-dimensional submanifold of 
$\mR^n$. Thus there is a smooth implicit definition of $E$, that is
to say a neighborhood $U$ of $E$ in $\mR^n$ and a smooth map 
$h : U \longrightarrow \mR^{n-2}$ of rank $n-2$ everywhere such that
$E=\ensemble{x\in U}{h(x)=0}$.

Moreover, assume there is a polynomial $P_{i_0}\notin\mR[P_1,P_2]$ (we 
choose it of minimal degree).

Let $d$ be the degree of $P_1$. We consider the taylor developpement of order
1 of $h$ in 0 and estimate it in 
$\left(P_1(x_1,x_2),\dots,P_n(x_1,x_2)\right)$. Noting $h_j$ the 
$j^{th}$ coordinate fonction of $h$ and $\partial_i$ the derivation in the
$i^{th}$ variable, we get for each $j$ :
$$
0=\sum_i \partial_i h_j(0) P_i(x_1,x_2) + o\left(\norme{x_1,x_2}^d\right)
$$
where the sum is taken over the $P_i$'s of degree $d$, hence
$$
\begin{pmatrix}\partial_1 h_1(0)    &\dots&\partial_n h_1(0) \\
               \vdots               &     &\vdots            \\
               \partial_1 h_{n-2}(0)&\dots&\partial_n h_{n-2}(0)
\end{pmatrix}
\begin{pmatrix} P_1 \\
                P_2 \mbox{ if it is of degree $d$, 0 otherwise}\\
                \vdots\\
                P_i\mbox{ if it is of degree $d$, 0 otherwise}\\
                \vdots\\
                P_n\mbox{ if it is of degree $d$, 0 otherwise}
\end{pmatrix}
=0.
$$
Each line in the second matrix is given by the coefficients of the 
polynomial.

First assume that $P_2$ and $P_{i_0}$ are both of degree $d$. Thus
the family of $P_i$'s of degree $d$ is of rank at least 3, hence the
jacobian matrix of $h$ at the point 0 is of rank at most $n-3$ which
prevent $h$ from being an implicit definition of $E$.

Next assume that $P_2$ is of degree $d$ and $P_{i_0}$ of degree
$d_0$ greater than $d$. Thus $h$ is of corank at least 2 at the point
 0: we have
two independent linear combinations of the $\partial_i h(0)$'s which must
be zero and involve only the indices $i$ of degree $d$ polynomials. But 
we can now use the Taylor developpement of order $d_0$ to get for each 
$j$:
$$ 0=\sum_i \partial_i h_j(0)P_i(x_1,x_2)+Q_j(x_1,x_2)$$
where the sum is taken over all polynomials of degree $d_0$ which are
not in $\mR[P_1,P_2]$ and $Q_j$ is a polynomial of degree $d_0$ of 
$\mR[P_1,P_2]$.
Let $S$ be, in the vector space of all homogenous 
polynomials of degree $d_0$, a supplementary of the space 
$\mR^{d_0}[P_1,P_2]$
of those of $\mR[P_1,P_2]$. Let $P'_i$ be the
projection of $P_i$ on $S$ along $\mR^{d_0}[P_1,P_2]$. Thus we have
for each $j$:
$$0=\sum_i \partial_i h_j(0)P'_i(x_1,x_2)$$
where the sum is taken over all polynomials of degree $d_0$ which are
not in $\mR[P_1,P_2]$. As before, it gives a linear combination of
the $\partial_i h(0)$'s which must be zero, and is independent of the two 
we get previously as $P_{i_0}\notin\mR[P_1,P_2]$. Hence $h$ is of 
corank at least 3 in O and the
contradiction holds as before.

We can use the same proof for the case when $P_2$ is of degree greater
than $d$.
\finpreuve

%%%%%%%%%%%%%%%%%%%%%%%%%%%%%%%%%%%%%%%%%%%%%%%%%%%%%%%%%%%%%%%%%%%%%%%%%
\subsection{Compactifications of the hyperbolic plane: the elliptic   
            case}
%%%%%%%%%%%%%%%%%%%%%%%%%%%%%%%%%%%%%%%%%%%%%%%%%%%%%%%%%%%%%%%%%%%%%%%%%

%%%%%%%%%%%%%%%%%%%%%%%%%%%%%%%%%%%%%%%%%%%%%%%%%%%%%%%%%%%%%%%%%%%%%%%%%
\subsubsection{Analytic non necessarily algebraic compactification}
\label{AnalComp}

	We shall start with a description of all analytic 
compactifications of 
$\mathcal{E}$ into a closed disk, in the following sense:

\begin{defi}
A differentiable \defini{compactification} of a differentiable action 
$\alpha$ of a Lie group $G$ on a manifold $M$ is a triple 
$(N,\phi,\overline{\alpha})$ where $N$ is a compact manifold 
with boundary, $\phi : M \longrightarrow N$ is an embedding and
$\overline{\alpha}$ is a differentiable action of $G$ on $N$
such that $\phi(M)$ is dense in $N$ and $\overline{\alpha}$
is a prolongation of the action induced by $\alpha$ on $\phi(M)$.
\end{defi}

The work of Schneider \cite{Schneider}, Stowe \cite{Stowe} exposed
by Mitsumatsu \cite{Mitsumatsu} gives immediately the classification
of all such compactifications, which we recall in what follows.

We shall use the following basis for \sldr:
$$ H=\begin{pmatrix} 1&0\\0&-1 \end{pmatrix},
   K=\begin{pmatrix} 0&-1\\1&0 \end{pmatrix},
   L=\begin{pmatrix} 0&1\\1&0 \end{pmatrix}.$$

The infinitesimal generators for the projective compactification
are given on $\mR\times\mR_+$ by
\begin{eqnarray}
\overline{K}_{1+} & = & 2\frac{\partial}{\partial x} 
                                                    \nonumber\\
\overline{H}_{1+} & = & 2\left((\sin x)(1+y)\frac{\partial}{\partial x}
  + (\cos x)(2y+y^2)\frac{\partial}{\partial y} \right)
                                                    \nonumber\\
\overline{L}_{1+} & = & 2\left((\cos x)(1+y)\frac{\partial}{\partial x}
  -(\sin x)(2y+y^2)\frac{\partial}{\partial y} \right).
                                                    \nonumber
\end{eqnarray}
and can be completed by adding a point at infinity.

\begin{theo}[\cite{Schneider}\cite{Stowe}\cite{Mitsumatsu}]\label{ClassAnal}
By pulling back the restriction of the vector fields 
$\overline{K}_{1+}$, $\overline{H}_{1+}$, $\overline{L}_{1+}$
to $\mR\times\mR\etoile_+$ by the map $F_n(x,y)=(x,y^n)$ 
where $n$ is a non-negative integer and by taking their continuous
prolongations,
we get analytic vector fields 
$\overline{K}_{n+}$, $\overline{H}_{n+}$, $\overline{L}_{n+}$
on $\mR\times\mR_+$. For any analytic compactifications
of $\mathcal{E}$ into a closed disc, there is an unique $n$ and a 
$\mR\times\mR_+$ chart in which
these vector fields are the infinitesimal generators of the 
compactified action. 
\end{theo}

For example, $\overline{K}_{2+}$, $\overline{H}_{2+}$, 
$\overline{L}_{2+}$ are the infinitesimal generators for
the conformal compactification.

%%%%%%%%%%%%%%%%%%%%%%%%%%%%%%%%%%%%%%%%%%%%%%%%%%%%%%%%%%%%%%%%%%%%%%%%%
\subsubsection{Analytic algebraic compactifications}

We shall now study the algebraic analytic compactifications of $\mathcal{E}$
into a closed disc, that is to say the elliptic orbits whose closure
is an analytic submanifold with boundary in the projective space
$\mP(V)$.

We prove a more precise version of the  theorem \ref{TheoEllip} exposed 
in the introduction:

\begin{theo}\label{TheoEllipDetail}
Let $O$ be the orbit of
$x = \left[u_q(\imath X+Y)^{\frac{n_q}{2}}
              (-\imath X+Y)^{\frac{n_q}{2}}\right]_q.$

If all the element of the family 
$(\frac{n_{q_+}-n_q}{2})_{q\in I(x)}$ are even, thus
$\adherence{O}$ is an analytic submanifold with boundary
and the action of \SLDR\ on this disk is
conjugate to the projective action.

If there exists some $q_{2+}$ in $I(x)$ such that 
$\frac{n_{q_+}-n_{q_{2+}}}{2}=1$, thus $\adherence{O}$ is an analytic 
submanifold with boundary and the action of \SLDR\ on 
this disk is conjugate to the conformal action.

In all the other cases, $\adherence{O}$ is not an analytic 
submanifold with boundary.
\end{theo}

\preuve
The methods used here will be useful through all the following
sections.

We shall first consider the case when all the numbers 
$\frac{n_{q_+}-n_q}{2}$, where $q$ is in $I(x)$,
are even.
A model for the projective compactification is given by
the closure in $\mP(\mR_2[X,Y])$ of the orbit of
$[X^2+Y^2]$, which is contained in the affine chart
$\ensemble{[aX^2+bXY+(1-a)Y^2]}{a,b\in \mR}$.
The map
\begin{eqnarray}
\varphi : \mP(\mR_2[X,Y]) & \longrightarrow & \mP(V) \nonumber\\
\left[aX^2+bXY+(1-a)Y^2\right]       & \longmapsto     & \left[u_q
  \left(a(1-a)-\frac{b^2}{4}\right)^{\frac{n_{q_+}-n_q}{4}}\right.
                                                     \nonumber\\
        & & \left.\phantom{\rlap{\ensuremath{
        \left(a(1-a)-\frac{b^2}{4}\right)^{\frac{n_{q_+}-n_q}{4}}}}}
        (aX^2+bXY+(1-a)Y^2)^{\frac{n_q}{2}}\right]_q \nonumber
\end{eqnarray}
is injective, analytic (thanks to the hypothesis)
 and realizes a conjugacy between the 
projective action and the dynamics on $\adherence{O}$.

Moreover, it is an immersion since, noting $s, t, u, v$ the coefficients
of the terms in 
$X^{n_{q_+}}, X^{n_{q_+}-1}Y, Y^{n_{q_+}}, XY^{n_{q_+}-1}$,
we have
$
\frac{\partial s}{\partial a}=\frac{n_{q_+}}{2}a^{\frac{n_{q_+}}{2}-1},
\frac{\partial u}{\partial a}=-\frac{n_{q_+}}{2}(1-a)^{\frac{n_{q_+}}{2}-1}
$
and
$
\frac{\partial s}{\partial b}=0,
\frac{\partial t}{\partial b}=\frac{n_{q_+}}{2}a^{\frac{n_{q_+}}{2}-1},
\frac{\partial u}{\partial b}=0,
\frac{\partial v}{\partial b}=\frac{n_{q_+}}{2}(1-a)^{\frac{n_{q_+}}{2}-1}.
$

Hence the differential of $\varphi$ is of rank 2 everywhere.

This proves that $\adherence{O}$ is an analytic submanifold with
boundary and at the same time that the action of
\SLDR\ on it is conjugate to the projective one.

Next we shall consider the case when there exists some $q_{2+}$ in $I(x)$ 
such that $\frac{n_{q_+}-n_{q_{2+}}}{2}=1$.
A model for the conformal action is given by the closure of
$\mH^+$ in the Riemmann sphere. We consider the map
\begin{eqnarray}
\psi : \adherence{\mH^+}      & \longrightarrow & \mP(V) \nonumber\\
       a+\imath b & \longmapsto     & \left[u_q
  b^{\frac{n_{q_+}-n_q}{2}}
  ((a+\imath b)X+Y)^{\frac{n_q}{2}}((a-\imath b)X+Y)^{\frac{n_q}{2}}
                                      \right]_q \nonumber
\end{eqnarray}
which is injective, analytic and realizes a conjugacy between the
conformal action and the dynamics on $\adherence{O}$. Notice that
$\psi(\infty)=[u_q X^{n_q}]_{q\in I_+(x)}$.

Moreover developping the expression of $\psi(a+\imath b)$, we see
that a coefficient is $n_{q_+}a$ and another is $u_{q_{2+}}b$, 
so $\psi$ is everywhere of rank 2 and we can conclude 
as before.

For the last case, we use Proposition \ref{raccourci}.
We denote by $\alpha$ the smallest odd element of the family 
$(\frac{n_{q_+}-n_q}{2})_q$, we denote by $q_{2+}$ an index realizing this minimum. 
By hypothesis $\alpha>1$. We can
write an element of $\adherence{O}$ under the form:
$ \left[u_q(\impart z)^{\frac{n_{q_+}-n_q}{2}}
  \left((\impart z^2+\repart z^2)X^2
        +2\repart z XY
        + Y^2\right)^{\frac{n_q}{2}}\right]_q.$
All coordinates are homogeneous polynomials in $x=\repart z$
and $y=\impart z$. Among them $P_1=x$ (we define it up to a multiplicative 
constant) is of minimal degree.
Among those which are not in $\mR[P_1]$, $P_2=y^2$ is of minimal degree.
But $P_3=y^{\alpha}\notin\mR[P_1,P_2]$ hence
$\adherence{O}$ is not a smooth submanifold of \mP(V), therefore not an 
analytic one.
\finpreuve
 
\begin{rema}\label{TheoEllipDetailRem}
In this proof
we can see more than stated: the embeddings $\varphi$
and $\psi$ extend respectively to embeddings of a projective plane 
(union of
the elliptic orbit of 
$x$, the hyperbolic orbit
of $\left[\left(-\frac{1}{4}\right)^{\frac{n_{q_+}-n_q}{4}}u_q
X^{\frac{n_q}{2}}Y^{\frac{n_q}{2}}\right]_q$ which is a
Moebius strip and their common border, the circular
orbit of $[u_q Y^{n_q}]_{q\in I_+(x)}$) and a
sphere (union of the elliptic orbits of $x$
and of
$\left[(-1)^{\frac{n_{q_+}-n_q}{2}}u_q(-X^2+Y^2)^{\frac{n_q}{2}}
\right]_q$ and of their common border,
the circular orbit of $[u_q Y^{n_q}]_{q\in I_+(x)}$).
\end{rema}

Moreover, we see that if we are in the third case, the map $\varphi$
is not analytic but is a \diffb{\frac{\alpha-1}{2}} embedding of
the projective action, so we can state the following fact
concerning the differentiable case for elliptic orbits:

\begin{theo}
The only algebraic differentiable compactifications of $\mathcal{E}$ 
are equi\-valent to the projective or to the conformal ones. In the 
projective case there exist \diffb{k} non-analytic realizations for
each finite $k$, but any \diffb{\infty} realization is in fact analytic.
In the conformal case any \diffb{1} realization is in fact analytic.
\end{theo}

%%%%%%%%%%%%%%%%%%%%%%%%%%%%%%%%%%%%%%%%%%%%%%%%%%%%%%%%%%%%%%%%%%%%%
\subsection{Hyperbolic case}
%%%%%%%%%%%%%%%%%%%%%%%%%%%%%%%%%%%%%%%%%%%%%%%%%%%%%%%%%%%%%%%%%%%%%

Here we shall consider the closure of a hyperbolic 2-dimensional orbit, 
which has the form
$$\adherence{O} = \ensemble{\left[u_q(t_1-t_2)^{\alpha_{q_+}-\alpha_q}
                      (t_1X+ Y)^{\alpha_q}(t_2X+ Y)^{n_q-\alpha_q}
                \right]_q}{t_1,t_2\in\mR\mP^1}.$$

\begin{theo}
If $O$ is a Moebius strip (\lat{i.e} for each $q$, $n_q$ is even, 
$\alpha_q=\frac{n_q}{2}$ and $\alpha_{q_+}-\alpha_q$ is even),
$\adherence{O}$ is an analytic submanifold; moreover its union
with the elliptic orbit of 
$\left[\left(-\frac{1}{4}\right)^{\frac{n_{q_+}-n_q}{4}}u_q
(X^{n_q}+Y^{n_q})\right]$ is still analytic and the dynamics
is conjugate to the projective action of \SLDR\ on the projective
plane.

If there is some $q_{2+}$ such that 
$\alpha_{q_+}-\alpha_{q_{2+}}=1$, $\adherence{O}$
is an analytic submanifold of $\mP(V)$ and its dynamics is
conjugate to the natural product action of \SLDR\ on
$\mR\mP^1\times\mR\mP^1$.

In all the other cases, $\adherence{O}$ is not an analytic
submanifold.
\end{theo}

\preuve
The first case is given by the map $\varphi$ of the previous section
(see Remark \ref{TheoEllipDetailRem}).

In the second case, we consider the map
\begin{eqnarray}
\psi : \mR\mP^1\times\mR\mP^1 & \longrightarrow & \mP(V) \nonumber\\
       (t_1,t_2)              & \longmapsto &
                 \left[u_q (t_1-t_2)^{\alpha_{q_+}-\alpha_q}
                 \left(t_1 X+Y\right)^{\alpha_q}
                 \left(t_2 X+Y\right)^{n_q-\alpha_q}\right]_q
                                                         \nonumber
\end{eqnarray}
which is analytic, injective as the orbit is by hypothesis a cylinder
and is an immersion as the coefficient of the terms
in $XY^{n_{q_+}}$ and $Y^{n_{q_{2+}}}$ of $\psi(t_1,t_2)$
are respectively
$\alpha_{q_+}t_1+(n_q-\alpha_{q_+})t_2$ and $t_1-t_2$, which gives a 
partial jacobian matrix
{\footnotesize $\begin{pmatrix} \alpha_{q_+} & n_q-\alpha_{q_+} \\
                                          1 & -1     
               \end{pmatrix}$} whose determinant is $-n_{q_+}\neq 0$.
Hence $\adherence{O}$ is an analytic submanifold (without boundary)
of $\mP(V)$ and (see the topological study) its dynamics is conjugate 
to the product action of \SLDR\ on $\mR\mP^1\times\mR\mP^1$.

For the last case we 
use Proposition \ref{raccourci}. 
The only polynomial of degree 1 among
the coordinates is $P_1=\alpha t_1 + \beta t_2$ where we write $\alpha$
for $\alpha_{q_+}$ and $\beta$ for $n_{q_+}-\alpha_{q_+}$. We can
next choose 
$P_2=\frac{\alpha(\alpha-1)}{2}t_1^2+\alpha\beta t_1 t_2 +
     \frac{\beta(\beta-1)}{2} t_2^2$.
Setting $P_2'=(t_1-t_2)^2$, an easy computation gives 
$\mR[P_1,P_2]=\mR[P_1,P_2']$.

If $\alpha=\beta$, as $\adherence{O}$ is assumed to be a cylinder there
must exist some index $q_0$ such that $\alpha_{q_+}-\alpha_{q_0}$ is odd.
Thus one of the coordinates has the form $(t_1-t_2)^{\alpha_{q_+}-\alpha_{q_0}}$
which is not in $\mR[P_1,P_2']$, hence from Proposition \ref{raccourci}
we conclude that $\adherence{O}$ is not an analytic submanifold of $\mP(V)$.

If $\alpha\neq\beta$, we see after an easy computation that the 
coordinate 
$P_3=\frac{\alpha(\alpha-1)(\alpha-2)}{6}t_1^3
     +\frac{\alpha(\alpha-1)}{2}\beta t_1^2 t_2
     +\alpha\frac{\beta(\beta-1)}{2} t_1 t_2^2
     +\frac{\beta(\beta-1)(\beta-2)}{6} t_2^3$
of the term $X^3Y^{n_q-3}$ is not in $\mR[P_1,P_2]$ and the conclusion
still holds.
\finpreuve

%%%%%%%%%%%%%%%%%%%%%%%%%%%%%%%%%%%%%%%%%%%%%%%%%%%%%%%%%%%%%%%%%%%%
\subsection{Parabolic case}
%%%%%%%%%%%%%%%%%%%%%%%%%%%%%%%%%%%%%%%%%%%%%%%%%%%%%%%%%%%%%%%%%%%%

We shall finally consider the closure of a parabolic orbit, which
has the form
$\adherence{O}=
  \ensemble{\left[u_q d^{n_q-n_{q_-}} (tX+Y)^{n_q}\right]_q}{
  d\in\overline{\mR}\mbox{ and } t\in\mR\mP^1}
$
where $d\in\overline{\mR}$ means $d$ is real or $\pm\infty$.

We shall prove some lemmas before stating the general result.
Let $q_{2-}$ (respectively $q_{2+}$) be an index such that 
$n_{q_{2-}}$ (respectively $q_{2+}$) is minimal (respectively
maximal) among
$n_q$'s greater than $n_{q_-}$ (respectively lesser than
$n_{q_+}$).

\begin{lemm}\label{lemm1}
If $n_{q_-}=0$ and $\adherence{O}$ is a smooth submanifold of $\mP(V)$,
we must have $n_{q_2-}=1$ and hence $\adherence{O}$ is a projective plane.
\end{lemm}

\preuve
We shall use Proposition \ref{raccourci} once again, around the point
$[u_q]_{q\in I_-}$ corresponding to $d=0$, $t=0$. The least-dimensional 
non-constant polynomial among the local coordinates is 
$P_1=d^{n_{q_{2-}}}$. There is no other polynomial of the same degree,
so we can choose $P_2=tP_1\notin\mR[P_1]$. If $n_{q_2-}>1$, one of
the coordinates can be written as $t^2P_1\notin\mR[P_1,P_2]$ 
and $\adherence{O}$
can not be a smooth submanifold of $\mP(V)$.
\finpreuve

\begin{lemm}\label{lemm2}
If $\adherence{O}$ is a smooth submanifold of $\mP(V)$, we must have
\begin{itemize}
  \item $n_{q_+}-n_{q_{2+}}=n_{q_{2-}}-n_{q_-}$,
  \item for each $q$, $n_{q_+}-n_{q_{2+}}$ divides $n_{q_+}-n_{q}$.
\end{itemize}
\end{lemm}

\preuve
We use Proposition \ref{raccourci} twice.

We first look around the point $[u_q Y^{n_q}]_{q\in I_-}$
to prove that for each $q$, $n_{q_{2-}}-n_{q_-}$ divides $n_q-n_{q_-}$. 
If $n_{q_-}=0$, we have $n_{q_{2-}}=1$ and the claim is obvious.
If $n_{q_-}>0$, we can choose $P_1=t$ and 
$P_2=d^{n_{q_{2-}}-n_{q_-}}$. For each $q$ there is a coordinate
which has the form $d^{n_q-n_{q_-}}$, hence by Proposition \ref{raccourci}
$n_{q_{2-}}-n_{q_-}$ must divide $n_q-n_{q_-}$.

In particular $n_{q_{2-}}-n_{q_-}$ divides $n_{q_+}-n_{q_{2+}}$.

We now look around the point $[u_q Y^{n_q}]_{q\in I_+}$, where
local coordinates are given by writting a point of
$\adherence{O}$ under the form 
$\left[u_q e^{n_{q_+}-n_q} (tX+Y)^{n_q}\right]_q$ after a change
of coordinates $e=d^{-1}$.
We can choose $P_1=t$ and $P_2=e^{n_{q_+}-n_{q_{2+}}}$, thus as
there is coordinates of the form $e^{n_{q_+}-n_q}$, for all $q$,
$n_{q_+}-n_{q_{2+}}$ divides $n_{q_+}-n_q$.

In particular $n_{q_+}-n_{q_{2+}}$ divides $n_{q_{2-}}-n_{q_-}$
and the conclusion holds.
\finpreuve

It is easy to see that the necessary conditions given in the previous
lemma are also sufficient if $n_{q_-}\neq0$ for $\adherence{O}$ to be 
an analytic submanifold of $\mP(V)$:
around each point of $\adherence{O}$ we can find local coordinates
of the form
$P_{k,l}=d^{k(n_{q_+}-n_{q_{2+}})}t^l$ where $k$ and $l$ are integers 
and for some coordinates we have $(k,l)=(0,1)$ or $(k,l)=(1,0)$, hence
writting $P_{k,l}-{P_{1,0}}^k {P_{0,1}}^l=0$ we get an analytic implicit
local definition of $\adherence{O}$.
If $n_{q_-}=0$ the combination of the conditions of the two lemmas are
also sufficient for $\adherence{O}$ to be analytic since we can find
local coordinates of the previous form or, around the points given by
$d=0$, of the form $P_{k,l}=d^k t^l$ with $k>0, k\geqslant l$; for some
coordinates we have $(k,l)=(1,1)$ and $(k,l)=(1,0)$ hence we get an
analytic implicit local definition of the form 
$P_{k,l}-{P_{1,1}}^l {P_{1,0}}^{k-l}=0$.

Moreover, if we map a point given by parameters $d,t$ from the closure
of an analytic parabolic orbit to the point given by the same parameters
on another such orbit closure of the same topology (projective plane, 
Klein bottle or cylinder) and with the same value for $n_{q_{2-}}-n_{q_-}$
we build an analytic diffeomorphism between them:
$$\left[u_q d^{n_q-n_{q_-}} (tX+Y)^{n_q}\right]_q 
  \longmapsto (d^{n_{q_{2-}}-n_{q_-}},t)
  \longmapsto \left[u'_q d^{n_q-n_{q'_-}} (tX+Y)^{n_q}\right]_q.$$

Finally, if we consider the differential in the point 
$x=[u_q Y^{n_{q_-}}]_{q\in I_-}$ of an element 
{\footnotesize$\begin{pmatrix} a&0\\c&a^{-1} \end{pmatrix}$} of the 
stabilizer of $x$ we find that its eigenvalues are $a^{-2}$ and
$a^{-(n_{q_{2-}}-n_{q_-})}$, so two closures of orbits with different
values of $n_{q_{2-}}-n_{q_-}$ can not be differentiably conjugate.
Hence we can state:

\begin{theo}
The conditions of Lemmas \ref{lemm1} and \ref{lemm2} are sufficient
for $\adherence{O}$ to be an analytic submanifold of $\mP(V)$.
Two analytic parabolic orbits are analyticaly conjugate if and only
if they have the same topology and the same value for
$n_{q_{2-}}-n_{q_-}$ (and they are not even differentiably conjugate 
otherwise). In particular there is one parabolic algebraic action on the
projective plane, a countable family of actions on the Klein bottle
and a countable family of actions on the closed cylinder.
\end{theo}

The last point we have to study in order to complete the proof of the 
results stated in the introduction is the way the cylindric orbits
are glued together.

Let $O$ be a cylindric analytic orbit associated to a projective
element $[u_q]_q$. Its boundary is the union of the two circular
orbits associated to the projective elements $[u_q]_{q\in I_+}$ and
$[u_q]_{q\in I_-}$, which we call respectively the \defini{upper
component} and the \defini{lower component} of the boundary.

An element of $\adherence{O}$ can be writen 
$\left[u_q d^{n_q-n_{q_-}} (tX+Y)^{n_q}\right]_q$ around the lower
component of the boundary. For each $q$ we denote by $k_q$ the integer 
$\frac{n_q-n_{q_-}}{n_{q_{2-}}-n_{q_-}}$.
The coordinates 
$c_{q,l}=u_q d^{n_q-n_{q_-}} t^l$ satisfy the implicit 
definition given previously:
$$\frac{1}{u_q}c_{q,l}-\frac{1}{u_{q_{2-}}}{c_{q_{2-},0}}^{k_q}
\frac{1}{u_{q_-} n_{q_-}} {c_{q_-,1}}^l=0.$$
Let $O'$ be the cylindric analytic orbit associated with the 
projective element $[u'_q]_q$ where $u'_q=(-1)^{k_q}u_q$.
Thus the lower component of its boundary is the same than for 
$O$ and as around it the coordinates of $O'$ satisfy the same
implicit parametrization, $O$ and $O'$ are analytically glued
together around their lower component.

With the same method we see that $O$ and the orbit $O''$ associated
with $[u''_q]_q$ where $u''_q=(-1)^{k_{q_+}-k_q}u_q$ are analytically
glued around their common upper component.

If $k_{q_+}$ is even $O'=O''$ and $O$ together with $O'$ gives a torus
with two open orbits, if $k_{q_+}$ is odd $O'\neq O''$ but they are
both glued analytically with $O'''$, the parabolic orbit associated with
$[(-1)^{k_{q_+}}u_q]_q$. Hence we have proven the last remaining result:

\begin{theo}
Let $O$ be a parabolic, cylindric, analytic orbit associated to
$[u_q]_q$.

If $k_{q_+}=\frac{n_{q_+}-n_{q_-}}{n_{q_{2-}}-n_{q_-}}$ is even,
the union of the two parabolic orbits associated to $[u_q]_q$
and $[(-1)^{k_q}u_q]_q$ is a torus analytically embedded in $\mP(V)$.

If $k_{q_+}$ is odd, the union of the four parabolic orbits associated
to $[u_q]_q$, $[(-1)^{k_q}u_q]_q$, $[(-1)^{k_{q_+}-k_q}u_q]_q$ and
$[(-1)^{k_{q_+}}u_q]_q$ is a torus analytically embedded in $\mP(V)$.
\end{theo}

\nocite{*}
\bibliographystyle{plain}
\bibliography{biblio}

\begin{thebibliography}{1}

\bibitem{Mitsumatsu}
Yoshihiko Mitsumatsu.
\newblock {${\rm SL}(2;{\bf R})$}-actions on surfaces.
\newblock In {\em Geometric study of foliations (Tokyo, 1993)}, pages 375--389.
  World Sci. Publishing, River Edge, NJ, 1994.

\bibitem{Schneider}
C.~R. Schneider.
\newblock {${\rm SL}(2,\,R)$} actions on surfaces.
\newblock {\em Amer. J. Math.}, 96:511--528, 1974.

\bibitem{Serre}
Jean-Pierre Serre.
\newblock {\em Alg\`ebres de {L}ie semi-simples complexes}.
\newblock W. A. Benjamin, inc., New York-Amsterdam, 1966.

\bibitem{Stowe}
Dennis~C. Stowe.
\newblock Real analytic actions of {${\rm SL}(2,{\bf R})$} on a surface.
\newblock {\em Ergodic Theory Dynam. Systems}, 3(3):447--499, 1983.

\end{thebibliography}
\addcontentsline{toc}{section}{Bibliographie}

\signature

%%%%%%%%%%%%%%%%%%%%%%%%%%%%%%%%%%%%%%%%%%%%%%%%%%%%%%%%%%%%%%%%%%%%%
\end{document}